\documentclass[conference,a4paper]{IEEEtran}

\usepackage{verbatim}
\usepackage{amsfonts,amsmath,mathrsfs,amssymb,amsbsy}
\usepackage[final]{graphicx}
\usepackage{times,cite}
\usepackage{balance}

\long\def\symbolfootnote[#1]#2{\begingroup%
\def\thefootnote{\fnsymbol{footnote}}\footnote[#1]{#2}\endgroup}

\usepackage{verbatim}
\usepackage[]{float,latexsym}
\usepackage{amsfonts,amsmath,mathrsfs,amssymb,amsbsy}
\usepackage{url}

\newtheorem{theorem}{Theorem}[section]

\newcommand{\Prob}{\mathsf{P}}
\newcommand{\Expect}{\mathsf{E}}
\newcommand{\indic}{\mathbb{I}}

\usepackage{color}
\definecolor{lightblue}{rgb}{.7, .8, 1}
\definecolor{lightgreen}{rgb}{.6, 1, .6}
\usepackage{color}
\definecolor{brown}{rgb}{1,0.38,0.03}

\definecolor{OliveGreen}{rgb}{.2,0.6,0.2}
\definecolor{BrickRed}{rgb}{.7,0.2,0.2}

\newcommand{\ignore}[1]{} 

\long\def\symbolfootnote[#1]#2{\begingroup%
\def\thefootnote{\fnsymbol{footnote}}\footnote[#1]{#2}\endgroup}

\DeclareMathOperator*{\argmax}{arg\,max}

\newcommand{\tauv}{\tau_{\scriptscriptstyle \text{V}}}
\newcommand{\tauu}{\tau_{\scriptscriptstyle \text{U}}}
\newcommand{\tauhb}{\tau_{\scriptscriptstyle \text{HB}}}

\begin{document}

\sloppy

\title{Quickest Hub Discovery in Correlation Graphs}

\author{
  \IEEEauthorblockN{Taposh Banerjee}\\
  \IEEEauthorblockA{School of Engineering and Applied Sciences\\
    Harvard University\\
    Cambridge, MA, USA\\
    Email: {tbanerjee}@seas.harvard.edu}
  \and
\IEEEauthorblockN{Alfred O. Hero III}\\
 \IEEEauthorblockA{Department of EECS\\
    University of Michigan\\
    Ann Arbor, MI, USA\\
    Email: {hero}@umich.edu}
}


\maketitle

\begin{abstract}
A sequential test is proposed for detection and isolation of hubs in a correlation graph. Hubs in a correlation graph of a random vector are variables (nodes) that have a strong correlation edge. It is assumed that the random vectors are high-dimensional and are multivariate Gaussian distributed. The test employs a family of novel local and global summary statistics generated from small samples of the random vectors. Delay and false alarm analysis of the test is obtained and numerical results are provided to show that the test is consistent in identifying hubs, as the false alarm rate goes to zero. 

%
\end{abstract}

\section{Introduction} \label{sec:Intro}
In this paper we consider the problem of quickest detection of a change in correlation between 
variables in a sequence of high-dimensional vectors. We are also interested in isolating or identifying 
the variables whose magnitude correlation with other variables has changed the most. We refer to such variables as \textit{hubs}. 
A precise definition will be given below. 

The problem of correlation detection and 
estimation is an important one in science and engineering. The problem 
is particularly difficult when the data involved is high-dimensional; see \cite{pourahmadi2013high}, 
\cite{fan2016overview}, and \cite{cai2016estimating} 
for a survey of the literature. We are interested in the sequential version of the problem where the variables are initially uncorrelated or independent, 
and due to some event, the correlation between the variables abruptly changes. This has applications  
in problems of fault detection, anomaly detection, detection of changes in time-series data or financial data, etc. 
For simplicity, we restrict our attention to multivariate Gaussian random vectors. However, the results in this paper 
are application to the broader class of elliptically contoured distributions \cite{bane-firo-hero-isit15}. 

If the pre-change parameter are known (means and variances of the uncorrelated variables), 
the classical approach to sequential detection of a change in the covariance matrix of a sequence of Gaussian random vectors 
is to use the generalized likelihood ratio (GLR) based Cumulative Sum (CUSUM) method of Lorden \cite{lord-amstat-1971}. 
In this algorithm one sequentially estimates the covariance of the data and applies classical quickest change detection (QCD) 
test using the estimated covariance. 
But, this algorithm works only for low dimensional data. One approach to handle high-dimensional data, 
is to make assumptions about the structure of the covariance of the matrix, like sparsity, and then estimate the high-dimensional covariance matrix. 
Although this approach is appealing, it is generally hard to design sequential tests based on such an approach. 
When the pre-change parameters are not known, as is the case in this paper, the problem is even harder. For low dimensional data, one can use a double-GLR approach \cite{lai-jrss-1995}, 
but optimality of such tests is still under investigation. 

In this paper we take a random matrix based approach to solving this problem. We collect the sequence of 
high-dimensional vectors in batches 
of successive samples from which we form a random matrix.
We then define summary statistics for the 
random matrix and obtain approximations for the asymptotic distribution for the summary statistics as the dimension 
of the vector goes to infinity. The asymptotic distribution we obtain belongs to a one-parameter exponential family, and
a change in correlation of the random vectors induces a change in the parameter of the random matrix. 
We use this fact to detect the change in correlation. This works differs from our effort in \cite{bane-firo-hero-isit15}
in that we also define a family of local summary statistics, one for each variable, and obtain their approximate asymptotic 
distributions. We use these local summary statistics to isolate or identify variables that have seen the most amount of change
in correlation with other variables. This fault isolation process is what we call hub discovery. 
We obtain asymptotic expressions for delay and false alarm  
for our proposed procedure, and show via numerical results that the test is consistent in hub discovery. 

\section{Problem Description}\label{sec:Prob}
We have a sequence of samples of independent and identically distributed (i.i.d.)
Gaussian $p$-variate random vectors $\{\mathbf{X}(k)\}$ with densities
\begin{equation}\label{eq:ellipticalDef}
f_{\mathbf{X}(k)}(\mathbf{x}) = \frac{1}{ \sqrt{2\pi}^p |\Sigma_k|^{p/2}} e^{\left(-\frac{1}{2} (\mathbf{x}-\boldsymbol{\mu}_k)^T\boldsymbol{\Sigma}_k^{-1}(\mathbf{x}-\boldsymbol{\mu}_k) \right)}. 
\end{equation}
For some time parameter $\gamma$ the samples are assumed to have common covariance 
parameter $ \mathbf{\Sigma}_0=D_p$ for $k < \gamma$, where $D_p$ is a digonal matrix,   
and common covariance parameter $\mathbf{\Sigma} \neq D_p$ and for $k\geq \gamma$. 
$\gamma$ is called the change point and the pre-change and post-change distributions of $\mathbf{X}(k)$ 
are denoted by $f_{\mathbf{X}}^0$ and  $f_{\mathbf{X}}^1$, respectively. 
No assumptions are made about the mean parameter $\boldsymbol{\mu}_k$, and can take different 
values for different $k$. 
More specifically, the change-point model is described by: 
\begin{equation}\label{eq:PreChg_PstChg_Den}
\begin{split}
\mathbf{X}(k) &\sim f_{\mathbf{X}}^{0}(\mathbf{x}) = \frac{1}{ \sqrt{2\pi}^p |\Sigma_0|^{p/2}} e^{-\frac{1}{2}(\mathbf{x}-\boldsymbol{\mu}_k)^T\boldsymbol{\Sigma}_0^{-1}(\mathbf{x}-\boldsymbol{\mu}_k)}, \\
&\; \hspace{7cm} k < \gamma, \\ 
       &\sim f_{\mathbf{X}}^{1}(\mathbf{x}) = \frac{1}{ \sqrt{2\pi}^p |\Sigma|^{p/2}} e^{-\frac{1}{2}(\mathbf{x}-\boldsymbol{\mu}_k)^T\boldsymbol{\Sigma}^{-1}(\mathbf{x}-\boldsymbol{\mu}_k)}, \\
       &\hspace{7cm} k \geq \gamma
\end{split}
\end{equation}
Thus, in the above change point model, the variables are \textit{uncorrelated before change}, and \textit{correlated after change}. 
The objective is to detect this change in correlation as quickly as possible, while avoiding false alarms, 
and also correctly isolating the variables that have experienced the most change in their correlation level (to be made precise below).  The mean parameters $\{\boldsymbol{\mu}_k\}$, the diagonal entries in $\boldsymbol{\Sigma}_0$ and the post-change
covariance matrix $\boldsymbol{\Sigma}$, are \textit{unknown}.

A decision-maker sequentially acquires the samples $\{\mathbf{X}(k)\}$.
At each time point $k$ the decision-maker decides to either stop sampling, 
declaring that the change has occurred, i.e., $k \geq \gamma$, or to continue sampling. 
The decision to stop at time $k$ is only a function of $(\mathbf{X}(1), \cdots, \mathbf{X}(k))$. 
Thus, the time at which the decision-maker decides to stop sampling, say $\tau$, is a stopping time 
for the vector sequence $\{\mathbf{X}(k)\}$. 
At the stopping time $\tau$, the decision maker incurs a delay of $\tau-\gamma$, if $\tau \geq \gamma$. 
If $\tau < \gamma$, we say that we have an event of false alarm. 
At time $\tau$, using the available data $\left(\mathbf{X}(1), \cdots, \mathbf{X}(\tau)\right)$, the decision maker 
also has to identify or isolate variables that have undergone the most amount of change in correlation.   
We refer to such variables as hubs. We now formalize the notion of a hub. 


For a non-diagonal covariance matrix $\boldsymbol{\Sigma}$ with correlation coefficients $\{\rho_{ki}\}$ define
\begin{equation}\label{eq:V_Sigma}
V_k(\boldsymbol{\Sigma}) = \max_{i \neq k} |\rho_{ki}|, \quad \mbox{for}\; k \in \{1, \cdots, p\},
\end{equation}
to be the maximum magnitude correlation coefficient for the $k$th variable. 
Hubs are defined as
\begin{equation}\label{eq:Hub_Hc}
\begin{split}
\mathsf{H} &=\{k: V_k = \max_{1 \leq j \leq p} V_j \}.
\end{split}
\end{equation}
Define the correlation graph $\mathcal{G}(\boldsymbol{\Sigma})$ 
associated with the post-change correlation matrix $\boldsymbol{\Sigma}$ 
as an undirected graph with $p$ vertices, each representing a variable in the vector $\mathbf{X}$. 
An edge is present between vertices $i$ and $j$ if the magnitude of the correlation coefficient $\rho_{ij}$
between the $i^{th}$ and $j^{th}$ components of the random vector $\mathbf{X}$ 
is nonzero, i.e., if $|\rho_{ij}| > 0$, $i \neq j$. The correlation graph $\mathcal{G}(\boldsymbol{\Sigma})$ can be treated as a weighted 
graph with weight on edge between node $i$ and $j$ equal to $|\rho_{ij}|$, provided the latter is nonzero. 
Then, hubs are the nodes with highest maximum weights. 

There are also other interesting ways to define 
hubs. For example, one can select a threshold $\rho$ and define hubs as those variables for which $V_k > \rho$. 
Another way to define hubs is to redefine correlation graphs with edges present only of $|\rho_{ij}| \geq \rho$, and 
then define hubs as the nodes with the highest degree. Results on hub discovery for these other type of hubs will be reported elsewhere. In this paper we only discuss hub discovery for hubs defined in \eqref{eq:Hub_Hc}.  

At time $\tau$, let 
\begin{equation}
\mathcal{D}\left(\mathbf{X}(1), \cdots, \mathbf{X}(\tau)\right) \in 2^{\{1,\cdots,p\}}
\end{equation}
 be the decision function 
that selects a subset of the $p$ variables as the hubs. 
The decision-maker's objective is to choose a test, i.e. $\{\tau, \mathcal{D}\}$, to 
detect this change in correlation of the random vectors, as quickly as possible, 
subject to a constraint on the rate of false alarms, and on the rate of false isolation. 

A problem formulation that captures the above performance trade-offs is the following. 
We seek $\{\tau, \mathcal{D}\}$ to solve 
\begin{equation}\label{prob:Pollak}
\begin{split}
 \min_{\{\tau, \mathcal{D}\}}   &  \quad  \sup_{\gamma\geq 1}  \; \Expect_\gamma[\tau-\gamma| \tau \geq \gamma] \\
 \mbox{subj. to}  &  \quad \Expect_\infty[\tau] \geq \beta, \\
 \mbox{and }  &  \quad  \sup_\gamma\Prob_\gamma [\mathcal{D} \neq \mathsf{H} ]\leq \zeta,
\end{split}
\end{equation}
where $\Expect_\gamma$ is the expectation with respect to the probability measure under which the change occurs at $\gamma$, 
$\Expect_\infty$ is the expectation when the change never occurs, and 
$\beta$ and $\zeta$ are user-specified constraints on the mean time to false alarm and the probability 
of false isolation, respectively. 

This problem formulation is also studied in \cite{tart-sqa-2008} in the context of detecting and isolating 
a single affected stream out of a finite number of independent streams of observations, when the exact
pre- and post-change distributions are known. 
Even with independent streams, and known distributions, the author in \cite{tart-sqa-2008} obtained a solution for only 
a relaxed version of the above problem, by replacing the false isolation probability metric by $\Prob_1 [\mathcal{D} \neq \mathsf{H} ]$. 

In our case there are three major challenges towards developing a test or extending/applying the results from \cite{tart-sqa-2008}:
\begin{enumerate}
\item The means $\{\mu_k\}$, covariance matrices $\boldsymbol{\Sigma}_0$ and $\boldsymbol{\Sigma}$ are unknown
\item The high dimension $p$ that is much larger than the total number of samples we will possibly accumulate. 
\item The streams (one for each dimension) in our observation models are not independent: In fact, we are trying to detect a change 
in the level of dependence between the streams. 
\end{enumerate}

In the sequential analysis literature, for such intractable problems, 
an alternative approach often taken is to propose a sub-optimal test, 
and provide a detailed performance analysis of the proposed test. The performance analysis is often used to 
design the test (to choose thresholds). The effectiveness of the proposed test is then verified by simulations. 
See for example \cite{sieg-venk-astat-1995}, \cite{desobry2005online}, and \cite{chen2016sequential}. 

In \cite{bane-firo-hero-isit15}, we proposed a test for the non-parametric family of vector elliptically contoured densities, of which the Gaussian case considered here is a special case, with the objective of minimizing 
detection delay while avoiding false alarms.  
In our test in \cite{bane-firo-hero-isit15}, we used a \textit{global} summary statistic to detect the change in correlation.
However, the global nature of our summary statistics in \cite{bane-firo-hero-isit15} does not allow for classification of where the changes occurred in the correlation matrix, i.e., correlation hub discovery.
In this paper we propose a family of $p$ \textit{local} summary statistics, one for each variable, 
and obtain asymptotic densities for these local statistics. We use these and the global 
statistic from \cite{bane-firo-hero-isit15} to perform sequential classification of the hubs in addition to detection of the existence of hubs.

\section{A Random Matrix Solution Approach}\label{sec:RandMatModel}
We propose a random matrix based solution to our detection problem. 
The decision maker analyzes the $p$-dimensional vector sequence $\{\mathbf{X}(k)\}$ in batches of size $n$, where $p \gg n$.  
This leads to a $n \times p$ random matrix sequence $\{\mathbb{X}(m)\}_{m \geq 1}$. 
Specifically, 
\begin{equation}\label{eq:mathbbXm}
\begin{split}
\mathbb{X}(1) &= [\mathbf{X}(1), \cdots, \mathbf{X}(n)]^T \\
\mathbb{X}(2) &= [\mathbf{X}(n+1), \cdots, \mathbf{X}(2n)]^T, \mbox{ etc.}
\end{split}
\end{equation}
Thus, for each $m$, each of the $n$ rows of the random matrix $\mathbb{X}(m)$ is an i.i.d. samples 
of the multivariate Gaussian $p$-variate random vector $\mathbf{X}$ $=[X_1, \cdots, X_p]^T$ \eqref{eq:ellipticalDef}.

A change in the covariance of the random vector sequence $\{\mathbf{X}(k)\}$ as per the description in \eqref{eq:PreChg_PstChg_Den} will 
change the covariance of the random matrix sequence $\{\mathbb{X}(m)\}$. We propose to detect the change by using a stopping rule for 
the sequence $\{\mathbb{X}(m)\}$. 
Note that:
\begin{enumerate}
\item A stopping rule for the sequence $\{\mathbb{X}(m)\}$ is also a stopping rule for the sequence $\{\mathbf{X}(k)\}$, but it takes values in multiples of $n$ (with respect to index $k$). 
The precise choice of $n$ will be discussed below. 
\item The law of the random matrix $\mathbb{X}$ \textit{at} the change point may consist of random vectors $\mathbf{X}$ from both pre- and post-change distributions. We ignore this issue in the rest of the paper and assume that all the vectors in the matrix either have density $f_{\mathbf{X}}^{0}$ or have density $f_{\mathbf{X}}^{1}$ \eqref{eq:PreChg_PstChg_Den}. However, 
see \cite{moustakides2016sequentially} and \cite{bane-chen-garcia-veer-icassp-2014}, where such transient behavior is addressed.
\end{enumerate}

In random matrix theory (RMT), asymptotic distributions are obtained for various functions of the random matrix, e.g., 
a semicircle law for the empirical eigen value distribution of a Wigner matrix, etc \cite{anderson2010introduction}. Such results are in general valid when both the dimensions 
of the matrix are taken to infinity. For our random matrix $\mathbb{X}$, however, the batch size $n$ is fixed, and the dimension $p$ is large. This is called the \textit{purely high-dimensional regime} \cite{hero2016foundational}. Thus, traditional RMT theorems cannot be applied in this regime. In this paper, we obtain approximate asymptotic distribution 
for a proposed summary statistic (in fact, a family of them) for the random matrix $\mathbb{X}$ in this regime. We then use this approximate distribution(s) to detect a change in the law of the random matrix. We also propose techniques to isolate the hubs.

\section{Local Summary Statistics for the Random Data Matrix}\label{sec:SumStat}
In this section we define a family of local summary statistics to detect hubs \eqref{eq:Hub_Hc}.
The summary statistics are a function 
of the sample correlation matrix obtained from the random data matrix $\mathbb{X}$. 
We then obtain their asymptotic distribution in the 
purely high dimensional regime of  $p\rightarrow \infty$, $n$ fixed. The asymptotic 
distributions are members of 
a one-parameter exponential family to be specified below.
In Section~\ref{sec:HubDiscAlgo}, 
we will use these summary statistics and their asymptotic distributions for sequential hub discovery.

For random data matrix $\mathbb{X}$ 
we write
$$\mathbb{X} = [\mathbf{X}_1, \cdots, \mathbf{X}_p] = [\mathbf{X}^T_{(1)}, \cdots, \mathbf{X}^T_{(n)}]^T,$$
where $\mathbf{X}_i = [X_{1i}, \cdots, X_{ni}]^T$ is the $i^{th}$ column and $\mathbf{X}_{(i)} = [X_{i1}, \cdots, X_{ip}]$  
is the $i^{th}$ row. 
Define the sample covariance matrix as 
$$ \mathbf{S} = \frac{1}{n-1} \sum_{i=1}^n (\mathbf{X}_{(i)} - \bar{\mathbf{X}})^T (\mathbf{X}_{(i)} - \bar{\mathbf{X}}),$$
where $\bar{\mathbf{X}}$ is the sample mean of the $n$ rows of $\mathbb{X}$.
Also define the sample correlation matrix as
$$ \mathbf{R} = \mathbf{D_S}^{-1/2} \mathbf{S} \mathbf{D_S}^{-1/2},  $$
where $\mathbf{D_A}$ denotes the matrix obtained by zeroing out all but the diagonal elements of the matrix $\mathbf{A}$.  
Note that, under our assumption that the covariance matrix $\mathbf{\Sigma}$ of the rows of $\mathbb X$ is positive definite, $\mathbf{D_S}$ is invertible with probability one.
Thus $\mathbf{R}_{ij}$, the element in the $i^{th}$ row and the $j^{th}$ column of the matrix $\mathbf{R}$, is the 
sample correlation coefficient between the $i^{th}$ and $j^{th}$ columns of $\mathbb X$.   

Define the summary statistics (see also \eqref{eq:V_Sigma})
\begin{equation}\label{eq:Sum_Stat_V}
V_k(\mathbb{X}) \; := \; V_k(\mathbf{R})\; := \; \max_{i \neq k}  |\mathbf{R}_{ki} |, \quad k \in \{1, \cdots, p\}.
\end{equation}
This defines a set of $p$ local statistics $\{V_k\}_{k \leq p}$, one for each random variable in the random vector.
We will use these statistics for detection of correlation strength hubs. 

\begin{figure}[htb]
\center 
\includegraphics[width=9cm, height=6cm]{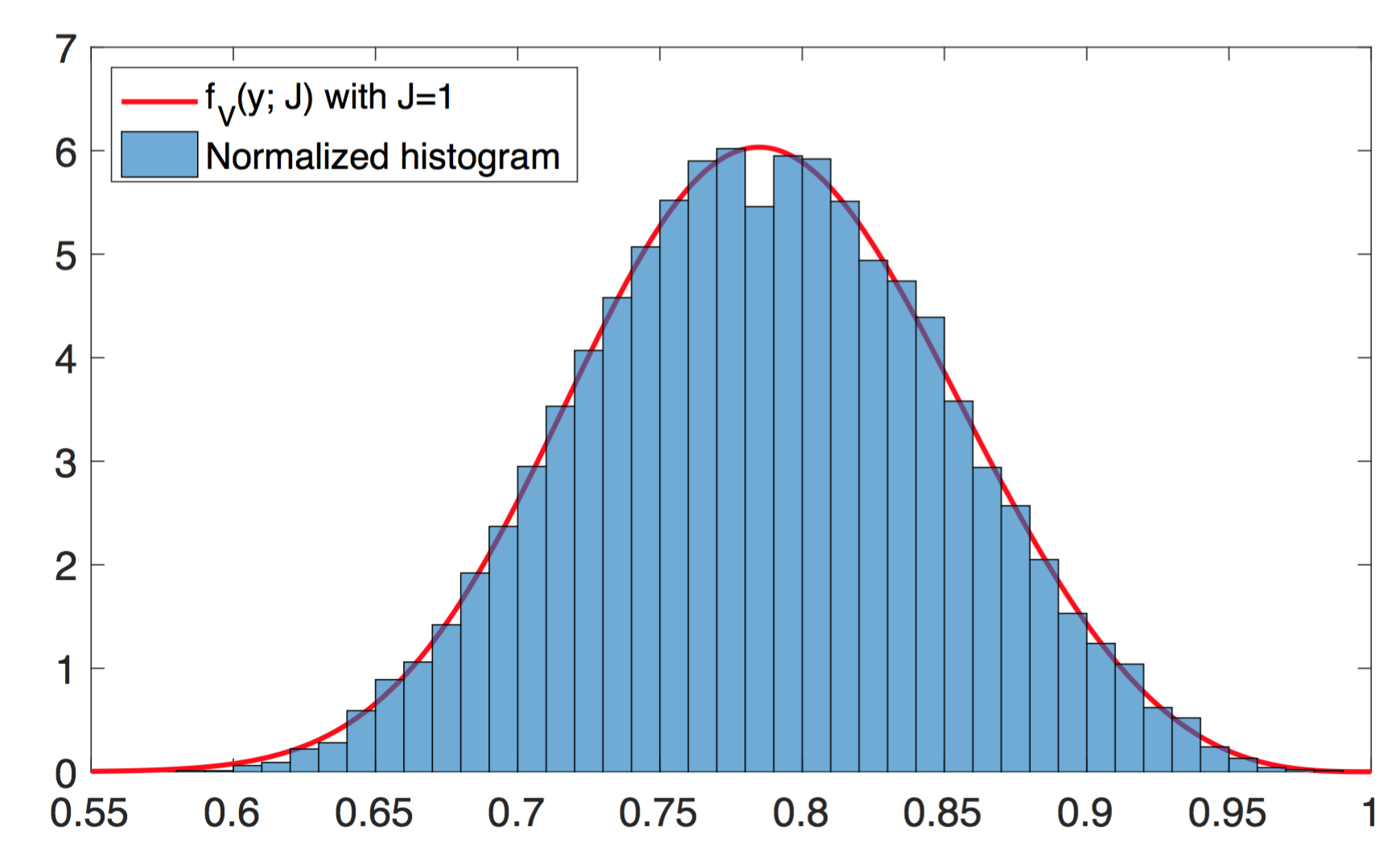}
\caption{Comparison of density $f_V$ in \eqref{eq:PDF_V_nok} for $J=1$ for multivariate Gaussian independent random variables with $n=10, p=100$ with the normalized histogram obtained from $10,000$ samples. Note that the population covariance is diagonal, yet the density is concentrated close to $1$.}
\label{fig:fV_histogram_Jeq1}
\end{figure}

Note that the summary statistics are the empirical estimates of the variables $\{V_k\}$ defined for the population correlation matrix $\boldsymbol{\Sigma}$ in \eqref{eq:V_Sigma}. If we could collect enough samples such that $n > p$, then we can rely 
on these estimates to check if they are close to zero, because the population version of these quantities are identically zero when the 
variables are uncorrelated. But, since $p \gg n$, we may have to collect a large number of samples before we can rely on these 
estimates. Such a delay may not be acceptable in practice, or those many samples may not even be available. 
In fact, for large $p$, $V_k(\mathbb{X})$ can be arbitrarily close to $1$; see \cite{hero-bala-IT-2012} and Fig.~\ref{fig:fV_histogram_Jeq1} in this paper. 

However, as we show below, one can obtain asymptotic law of these estimates as $p \to \infty$. The precise statement is provided 
in the theorems below. A change in correlation between 
the variables may be reflected in the asymptotic law of these estimates. This phenomena can be utilized for detecting 
a change in correlation. The rest of this section is devoted to obtaining the asymptotic laws. The algorithms are 
discussed in Section~\ref{sec:HubDiscAlgo}.  

We say that a matrix is row sparse of degree $j$ if there are no more than $j$ nonzero entries in any row. We say that 
a matrix is block sparse of degree $j$ if the matrix can be reduced to block diagonal form having a single $j \times j$ block, 
via row-column permutations. For a threshold parameter $\rho \in [0,1]$ define the correlation graph $\mathcal{G}_\rho(\boldsymbol{\Sigma})$ 
associated with the post-change correlation matrix $\boldsymbol{\Sigma}$ 
as an undirected graph with $p$ vertices, each representing a variable in the vector $\mathbf{X}$. 
An edge is present between vertices $i$ and $j$ if the magnitude of the correlation coefficient $\rho_{ij}$
between the $i^{th}$ and $j^{th}$ components of the random vector $\mathbf{X}$ 
is greater than $\rho$, i.e., if $|\rho_{ij}| \geq \rho$, $i \neq j$. 
We define $d_k$ to be the degree of vertex $k$ in the graph  $\mathcal{G}_\rho(\boldsymbol{\Sigma})$:
\begin{equation}\label{eq:d_Sigma}
\begin{split}
d_k(\boldsymbol{\Sigma}, \rho) \; &:= \; \# \{i \neq k:  |\mathbf{\rho}_{ki} | \geq \rho \}, \\
&\quad \quad  \quad k \in \{1, \cdots, p\}, \; \rho \in [0,1],
\end{split}
\end{equation}
where $\#A$ is the size of set $A$. 
Finally, throughout the index $k$ will be associated with the $k$th variable.

\medskip
\begin{theorem}\label{thm:CorrScr}
Let $\mathbf{\Sigma}$ be the covariance matrix of $\mathbf{X}$ or the rows of $\mathbb{X}$, $\mathbf{\Sigma}$ could be diagonal or non-diagonal depending on whether it is from the pre- or post-change data, respectively. Let $\mathbf{\Sigma}$ be row sparse of degree $j=o(p)$. Also let $p \to \infty$ and $\rho = \rho_p \to 1$ such that 
$p(1-\rho^2)^{n/2} \to 0$. 
Then,
\begin{equation}
\Prob(V_k(\mathbb{X}) \leq \rho ) \; \to \; e^{- \Lambda_k },  
\end{equation}
where 
$$\Lambda_k = \lim_{p\to \infty, \rho \to 1} \Lambda_{k,\rho} ,$$
with
$$\Lambda_{k,\rho} = (p-1)  P_0(\rho) J_k,$$
$$P_0(\rho) = I_{1-\rho^2} ((n-2)/2, 1/2),$$
with $I_x(a,b)$ being the regularized incomplete beta function with parameter $a$ and $b$.
The parameter $J_k$ is a positive real number that is a function of the asymptotic joint density of $\mathbf{X}$. 
Further, if the covariance matrix $\mathbf{\Sigma}$ of the p-variate vector $\mathbf{X}$ is block sparse of degree $j$, 
then 
$$J_k=1 + O((j/p)^{2}).$$ 
In particular, if $\mathbf{\Sigma}$ is diagonal then $$J_k=1.$$
\end{theorem}
\begin{IEEEproof}
The events $\{V_k(\mathbb{X}) \geq \rho\}$ and $\{d_k > 0\}$ are equivalent. Hence
\begin{equation}\label{eq:V_eq_N}
\Prob(V_k(\mathbb{X}) \geq \rho) = \Prob(d_k(\mathbb{X}, \rho) > 0). 
\end{equation}  
The result now follows from Theorem~\ref{thm:degree} below. 
\end{IEEEproof}
\medskip
\begin{theorem}[\cite{firouzi2013local}]\label{thm:degree}
Let $\mathbf{\Sigma}$ be row sparse of degree $j=o(p)$. Also let $p \to \infty$ and $\rho = \rho_p \to 1$ such that 
$p(1-\rho^2)^{n/2} \to 0$. Then
 $d_k(\mathbb X,\rho)$ converges in distribution to a Poisson random variable with rate parameter $\Lambda_k$. \end{theorem}
\medskip

Based on the asymptotic distribution obtained in Theorem~\ref{thm:CorrScr}, 
the  large $p$  distribution of $V_k$ defined in \eqref{eq:Sum_Stat_V} can be approximated by
\begin{equation}\label{eq:CDF_V}
\begin{split}
\Prob(V_k(\mathbb{X}) \leq \rho) &\approx \exp(-\Lambda_{k,\rho}),\\
&= \exp(-(p-1)  P_0(\rho) J_k) \; \rho \in [0,1]. 
\end{split}
\end{equation}
where $\Lambda_{k,\rho}$ is as defined in Theorem~\ref{thm:CorrScr}. Using \cite[Thm. 3.1]{firouzi2013local} it can easily be shown that the approximation error associated with (14) decays to zero at least as fast as  $p(1-\rho^2)^{n/2}\Delta$ where $\Delta$ is a dependency coefficient associated with the set of U-scores..    

The distribution \eqref{eq:CDF_V} is differentiable everywhere except at $\rho=0$ 
since $P(V_k(\mathbb{X})=0)>0$. For $\rho>0$ and large $p$,  $V_k$ has density  
\begin{equation}\label{eq:PDF_V_1}
f_{V_k}(\rho) \approx - (p-1) P_0'(\rho)  J_k \exp(-(p-1)  P_0(\rho) J_k), \; \rho \in (0,1]. 
\end{equation}
Note that $f_V$ in \eqref{eq:PDF_V_1} is the density of the Lebesgue continuous component of 
the distribution \eqref{eq:CDF_V} and that it integrates to $1-O(e^{-p^2})$ over $\rho\in (0,1]$.

For each $k$, the density $f_{V_k}$ is a member of a one-parameter exponential family with $J_k$ as the unknown parameter.  
Then, since $P_0(\rho)' = -\frac{2 (1-\rho^2)^{(n-4)/2}}{B((n-2)/2, 1/2)}$, where $B(a,b)$ is the beta function, we have for $ \rho \in [0,1]$,
the exponential family form of the density $f_V$ with parameter $J_k$:
\begin{equation}\label{eq:PDF_V}
\begin{split}
f_{V_k} &(\rho; J_k)  \\
=&  \; \frac{2 \; (p-1) \; J_k \; (1-\rho^2)^{\frac{n-4}{2}}}{B((n-2)/2, 1/2)}  \exp(-(p-1)  \; J_k   \; P_0(\rho) ).
\end{split}
\end{equation}

Note that the density $f_{V_k}$ corresponds to the $k$th variable, and the whole family $\{f_{V_k}\}$ depends on the index $k$ only through the parameter $J_k$, which can be different for different $k$. 
We thus rewrite the family in \eqref{eq:PDF_V} (with complete specification for easy accessibility) as
\medskip
\begin{equation}\label{eq:PDF_V_nok}
\begin{split}
f_{V} &(y; J)  =  J \; C(p, n) \;(1-y^2)^{\frac{n-4}{2}} e^{-(p-1)  \; J   \; P_0(y) }, 
\end{split}
\end{equation}
\medskip
where 
\begin{equation}
\begin{split}
C(p,n) &= \frac{2 \; (p-1)}{B((n-2)/2, 1/2)}, \\
P_0(y) &= I_{1-y^2} ((n-2)/2, 1/2),\\
I_x(a,b) & = \frac{B(x; a,b)}{B(a,b)},\\
B(x, a,b) & = \int_0^x t^{a-1} (1-t)^{b-1} dt,\\
B(a,b) & = B(1,a,b).
\end{split}
\end{equation}

\begin{figure}[htb]
\center 
\includegraphics[width=9cm, height=6cm]{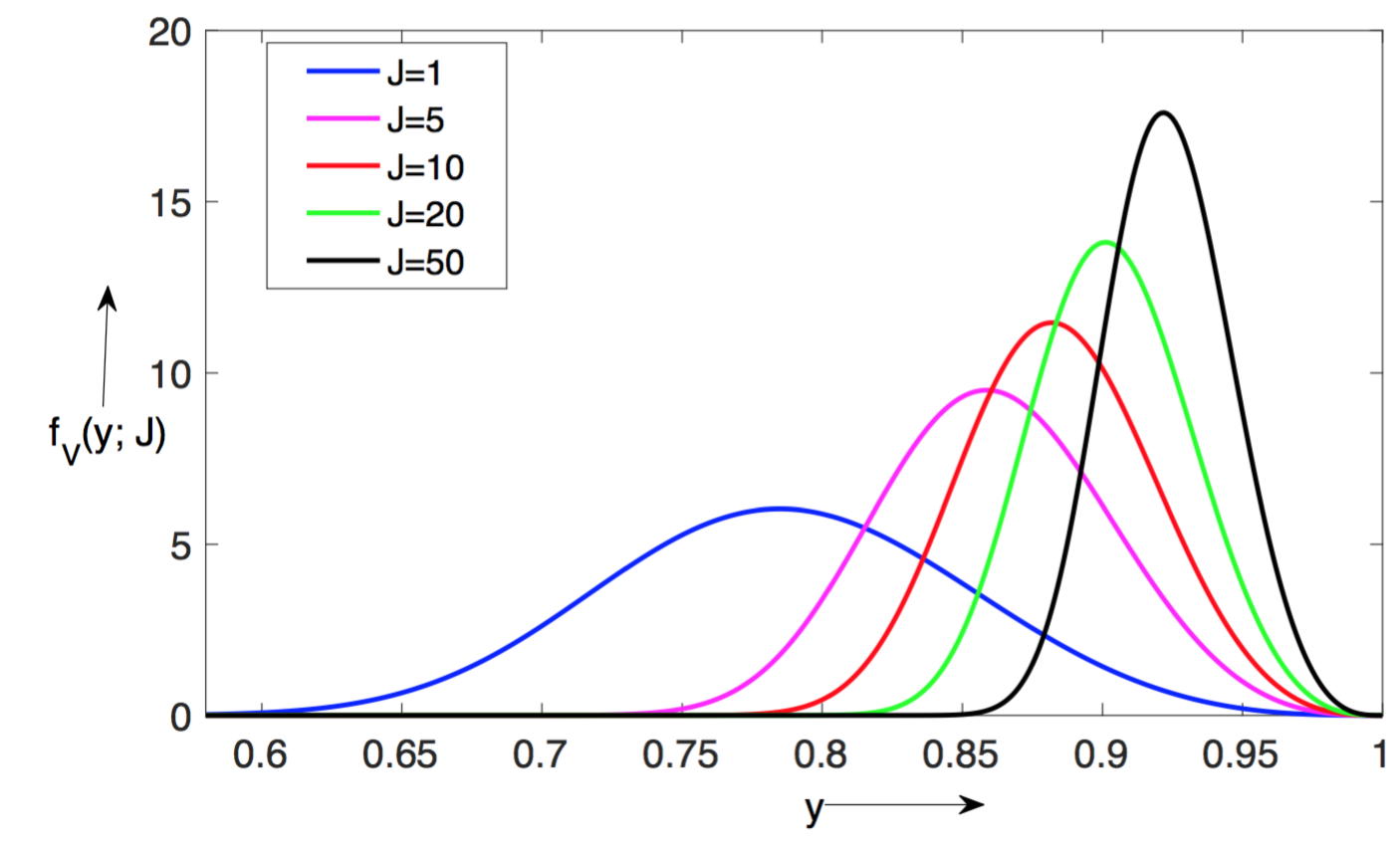}
\caption{Plot of density $f_V$ in \eqref{eq:PDF_V_nok} for various values of the parameter $J$ for $n=10, p=100$. }
\label{fig:fV_density}
\end{figure}

In Fig.~\ref{fig:fV_density} we have plotted the density $f_V$ for various values of $J$ for 
$n=10$, and $p=100$. We note that for the chosen values of $n$ and $p$, the density is concentrated close to $1$,
consistent with large values of $\rho$ arising in the purely high dimensional regime assumed in Theorem~\ref{thm:CorrScr} \cite{hero2016foundational}. Note that \textit{the variable $V$ tends to take on higher values as the  parameter $J$ increases}. 
This is the fundamental idea used for hub discovery in this paper. 

In Fig.~\ref{fig:fV_histogram_Jeq1} we have plotted the normalized histogram of the samples $V_1(\mathbb{X})$ for multivariate 
Gaussian independent random vectors and compare it with the approximation $f_V(\cdot\; ; 1)$. Note that the approximation is quite accurate. Also, note that for this case the population $V_1 = 0$, but the samples are concentrated in the interval $(0.55,0.95)$. This is due to the fact that $p \gg n$. 

\section{Global Summary Statistics}\label{sec:globalSumStat}
In \cite{bane-firo-hero-isit15} we proposed a global summary statistic and obtained its asymptotic distribution. 
Specifically, we showed that the statistic
\begin{equation}\label{eq:globalU}
U(\mathbb{X}) = \max_k V_k (\mathbb{X}),
\end{equation}
in the limit $p \to \infty$ (taken in a particular sense) and under sparsity assumptions similar to that in Theorem~\ref{thm:CorrScr}, has approximate density in an exponential family given by
\begin{equation}\label{eq:globalgU}
g(u;  \theta) = \frac{D \; \theta }{2} (1-u^2)^{\frac{n-4}{2}} \exp\left(-\frac{D}{2} \theta \;T(u) \right), \; u \in (0,1],
\end{equation}
where
\begin{equation}
\begin{split}
D=D_{p,n}\; &=\; 2p (p-1) B((n-2)/2, 1/2),\\
T(u) &= \int_{u}^1 (1-s^2)^{\frac{n-4}{2}} ds.
\end{split}
\end{equation}
The parameter $\theta$ is a positive parameter and has similar interpretation and nature as that of the parameter $J$ 
in \eqref{eq:PDF_V_nok}. Specifically, if the rows of $\mathbb{X}$ have diagonal covariance $\mathbf \Sigma$ then $\theta=1$, 
otherwise $\theta \neq 1$. 
Compared to $V_k (\mathbb{X})$, $U(\mathbb{X})$ has a global view of the covariance matrix. If any of the $V_k$s undergo
a change, this change will be reflected through a change in the distribution of $U$.  
However, it is not possible to isolate hubs based only on the statistic $U$. 
Sequential tests for detecting a change in correlation using the statistic $U$ is studied in detail in \cite{bane-firo-hero-isit15}. The sequential test from \cite{bane-firo-hero-isit15} will play a fundamental role in what follows. 

\section{Multivariate Quickest Hub Discovery}\label{sec:HubDiscAlgo}
We now discuss how to use the results from the previous sections to detect and isolate correlation hubs. 

Recall that we have a sequence of $p$-variate Gaussian distributed vectors $\{\mathbf{X}(k)\}$ with a change point model given in \eqref{eq:PreChg_PstChg_Den} and $p \gg 1$. Thus, the vectors are uncorrelated before change and correlated with covariance matrix $\mathbf{\Sigma}$ after change. We collect these vectors in batches of size $n$ and obtain a sequence of $n \times p$ random 
matrices $\{\mathbb{X}(m)\}$; see \eqref{eq:mathbbXm}. 
For quick detection, we should choose a batch size $n$. Thus, we also have $p \gg n$. 

Assuming that the post-change covariance matrix $\mathbf{\Sigma}$ is row sparse, a change in the law of 
$\mathbf{X}$ changes the law of the random matrix $\mathbb{X}$, which in turn affects the law of summary statistics
$U(\mathbb{X})$ and $\{V_k(\mathbb{X})\}_{k=1}^p$. Since, the distributions of these statistic each belong to a one-parameter 
exponential family \eqref{eq:PDF_V_nok} and \eqref{eq:globalgU}, a change in law within this family is a simple change 
in the value of the parameters $\{J_k\}$ and $\theta$. We propose to detect the change in correlation in $\mathbf{X}$ through a change 
in the parameters $\{J_k\}$ and $\theta$

From \eqref{eq:PDF_V_nok} and \eqref{eq:globalgU}, 
it follows that we can choose $n > 4$, because the densities are well defined for these values of $n$. 
Note also that here we assume that the change in correlation in $\mathbf{X}$ affects the parameters $\{J_k\}$ and $\theta$, otherwise our test is insensitive to the change in correlation.

\subsection{Correlation Change Detection Using Global Statistic}

Assuming that the post-change covariance matrix $\Sigma$ is row sparse, 
we map the random matrix sequence $\{\mathbb{X}(m)\}$ to the summary statistic sequence $\{U(\mathbb{X}(m))\}$.
The random variables $\{U(\mathbb{X}(m))\}$ have approximate 
density given by \eqref{eq:globalgU} which is in an exponential family with parameter $\theta$ which is $1$ before the 
change and some value $\theta \neq 1$ after change.
For simplicity, we refer to $U(\mathbb{X}(m))$ simply by $U(m)$. 
The QCD problem on the density $f_{\mathbf{X}}$, depicted in \eqref{eq:PreChg_PstChg_Den}, 
is reduced to the following QCD problem
\begin{equation}\label{eq:QCDProb_U}
\begin{split}
U(m) &\sim g(\cdot \; ; 1), \quad\hspace{1.25cm} m < \gamma \\ 
    &\sim g(\cdot \; ; \theta), \; \quad \theta \neq 1, \; m \geq \gamma. 
\end{split}
\end{equation}

Consider the following generalized likelihood ratio (GLR) based QCD tests (see \cite{lord-amstat-1971}, \cite{bane-firo-hero-isit15}) defined by the stopping time:
\begin{equation}\label{eq:tauU}
\begin{split}
&\tauu \\
&=\min \left\{m:\max_{1 \leq \ell \leq m} \sup_{\theta: |\theta-1| \geq \epsilon_u} \sum_{i=\ell}^m \log \frac{g(U(i); \theta)}{g(U(i); 1)} > A_u\right\},
\end{split}
\end{equation}
where we assign $\tauu=\infty$ if the right hand side of (23) is the empty set. Here $A_u > 0$ and $\epsilon_u>0$ are user-defined parameters, and $\epsilon_u$ represents the minimum magnitude of change, away from $\theta=1$, that the user wishes to detect. 

\subsection{Correlation Change Detection Using Local Statistics}
We can also detect a change in correlation using the local sample correlation statistics $\{V_k(\mathbb{X})\}_{k=1}^p$ defined in \eqref{eq:Sum_Stat_V}. 
For simplicity, $V_k(\mathbb{X}(m))$ is denoted simply by $V_k(m)$. 
For large $p$ and row sparse covariance, the random variables $\{V_k(m)\}$ have approximate 
density given by \eqref{eq:PDF_V_nok} which is in an exponential family with parameter $J_k$. This parameter is $1$ before the 
change and takes on some value $J_k \neq 1$ after the change. 
A natural locally pooled sqequential change detection test is to implement $p$ local sequential change detection tests in parallel, one for each variable $\{V_k(\mathbb{X})\}_{k=1}^p$ (equivalently, each variable in $\mathbf{X}$). An alarm is raised 
when \textit{any} of the local tests raises an alarm. 

Mathematically, the QCD problem on the density $f_{\mathbf{X}}$, depicted in \eqref{eq:PreChg_PstChg_Den}, 
is reduced to the following family of QCD problems (one for each $k$):   
\begin{equation}\label{eq:QCDProb_1toother}
\begin{split}
V_k(m) &\sim f_V(\cdot \; ; 1), \quad\hspace{1.25cm} m < \gamma \\ 
    &\sim f_V(\cdot \; ; J_k), \; \quad J_k \neq 1, \; m \geq \gamma.
\end{split}
\end{equation}
Consider the following GLR-based QCD tests defined by the family of local stopping times \cite{lord-amstat-1971}:
\begin{equation}\label{eq:taugk}
\begin{split}
&\tauv^{(k)} =\\
&\min  \left\{m: \max_{1 \leq \ell \leq m} \sup_{J_k: |J_k-1| \geq \epsilon_v} \sum_{i=\ell}^m \log \frac{f_{V}(V_k(i);J_k)}{f_{V}(V_k(i); 1)} > A_v\right\},
\end{split}
\end{equation}
where, again we assign $\tauv^{(k)}=\infty$ if the right hand side of \eqref{eq:taugk} is the empty set, and 
$A_v > 0$ and $\epsilon_v>0$ are user-defined parameters. 
This stopping time can be used to detect a change in the parameter $J_k$ for the $k$th variable. Define the overall local stopping time to be
\begin{equation}\label{eq:taug}
\tauv = \min_k \tauv^{(k)} = \min  \left\{m:  \max_{1 \leq k \leq p} G_k (m) > A_v\right\},
\end{equation}
where 
\begin{equation}\label{eq:GLRStat}
G_k (m) := \max_{1 \leq \ell \leq m} \sup_{J_k: |J_k-1| \geq \epsilon_v} \sum_{i=\ell}^m \log \frac{f_{V}(V_k(i);J_k)}{f_{V}(V_k(i); 1)}
\end{equation}
is the GLR statistic. 
The test $\tauv$ can be used to detect if the parameter of any of the $p$ variables is affected by the change. $\epsilon_v$ represents the minimum magnitude of change, away from $J_k=1$, that the user wishes to detect.

\subsection{Joint Detection and Hub Discovery}
For joint detection and hub discovery, we combine the local and global tests and propose to use the following maximum of the stopping times:
\begin{equation}\label{eq:tauHB}
\tauhb = \max\{ \tauv, \tauu\}. 
\end{equation}
Let $D_k(\mathbb{X}(1), \cdots, \mathbb{X}(\tauhb))$ be the binary decision variable with value $1$ if variable $k$ is declared a hub. We use the following rule for hub discovery: fix a positive integer $q$ and set
\medskip
\begin{equation}\label{eq:FaultIsoRule}
D_k(\mathbb{X}(1), \cdots, \mathbb{X}(\tauhb)) = \indic_{\left\{    G_k (\tauhb) \mbox{ is top $q$ statistic} \right\}}.
\end{equation}
\medskip

\begin{figure}[htb]
\center 
\includegraphics[width=8cm, height=5cm]{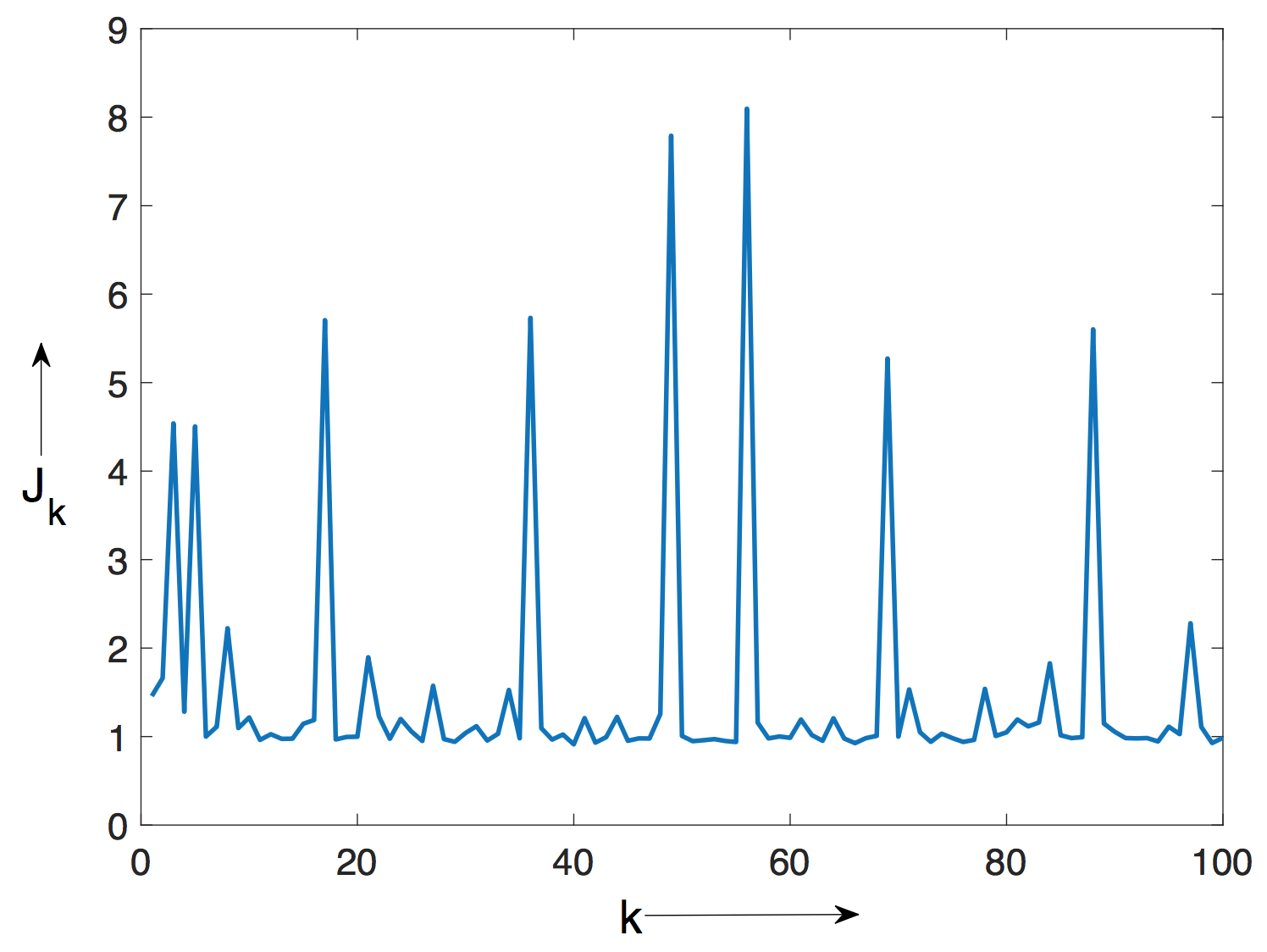}
\caption{Estimated $J_k$ values from a $p \times p$ row-sparse covariance matrix with $p=100$. Here $k$ is the index of the variables in $\mathbf{X}$. Note that the values of parameter $J_k$ for non-hub variables are also significant. Hubs are the two variables with the highest values of the parameter. }
\label{fig:JkPlot}
\medskip
\includegraphics[width=8.5cm, height=5.5cm]{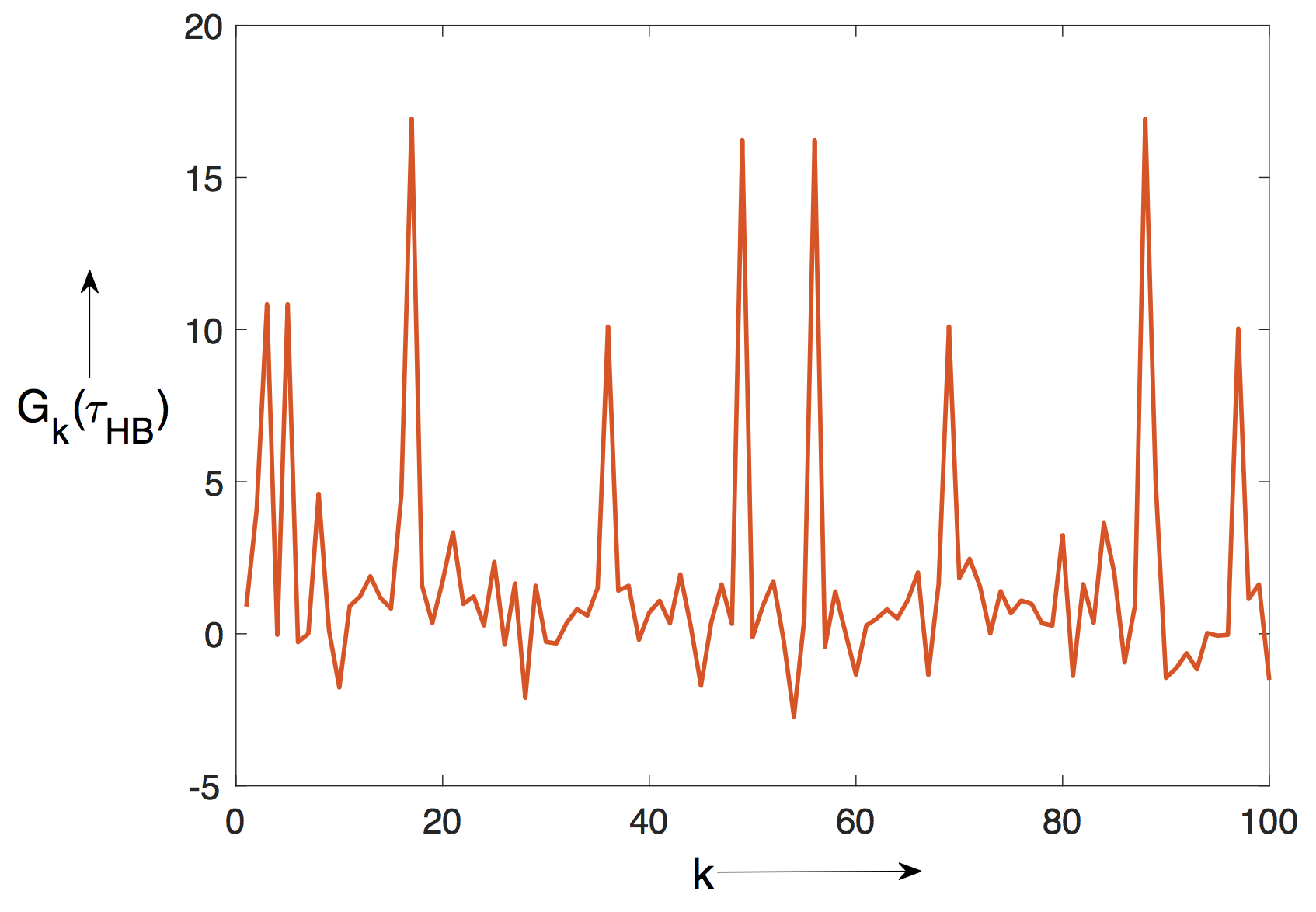}
\caption{Values of GLR statistic $G_k$ at the time of stopping $\tauhb$ corresponding 
to the model used in Fig.~\ref{fig:JkPlot}. Note that the $G_k$ values for the non-hub variables are comparable to those of the hub variables. In practice, one can just plot these $G_k$ values to capture all
the variables that are affected by a change in correlation. A close mathematical equivalent to this idea is the pick-the-top $q$ approach \eqref{eq:FaultIsoRule}.}
\label{fig:GkPlot}
\end{figure}

To appreciate the above rule, recall that a correlation strength hub is defined as the set of those variables that achieve the maximum in \eqref{eq:globalU}; see \eqref{eq:Hub_Hc}. Thus, there are at least two variables that are correlation strength hubs in the correlation graph. These hub variables in general have value of parameter $J_k$ that is higher than the parameter values for other variables.
However, there may be many variables that have experienced a change in correlation, and their $V_k$ value, and hence 
the value of parameter $J_k$, might be close to the values attained by the hubs. 
To illustrate this idea, we have plotted in Fig.~\ref{fig:JkPlot} estimated $\{J_k\}$ values from samples of a multivariate Gaussian vector with a row sparse covariance matrix (see Section~\ref{sec:NumericalResults}
on how the covariance matrix is generated). Here the hubs are the two variables with the highest values of $J_k$. 
But, as seen in the figure, there are many other variables that have experiences significant changes in correlation, and 
hence also have high values of $J_k$.  
In Fig.~\ref{fig:GkPlot}, we have plotted the values of the GLR statistic $\{G_k\}$ at the time of stopping $\tauhb$. 
As seen in the figure, in addition to the $G_k$ values of hubs, the values corresponding to other variables 
are also significant. Thus, a simple strategy like choosing the variable with the highest $G_k$ value as hub, will not work well. 

However, what is clear from Fig.~\ref{fig:JkPlot} and Fig.~\ref{fig:GkPlot} is that high $G_k$ values 
at the time of stopping correspond to high $J_k$ values of the summary statistics, and these values together have significantly higher magnitude than the values for the rest of the variables. This motivates
the pick-the-top approach. 
In fact, in practice, \textit{one can plot the $G_k$ values this way and identify all the variables that have experienced a change in correlation}. 

Picking the variables corresponding to the top $q$ values where $q > 2$ ensures that the hubs are discovered with high probability, 
provided a slightly modified definition of hub discovery is used: the variables in the hub $\mathsf{H}$ are considered discovered if
\begin{equation}\label{eq:realDefHubDiscovery}
\mathsf{H} \subset \{k : \mathcal{D}_k= 1\}.
\end{equation}
Indeed, this fact is used in Section~\ref{sec:NumericalResults} for numerical computations,
where it is shown that the proposed rule's probability of false hub discovery goes to zero, as the false alarm rate 
goes to zero; see Fig.~\ref{fig:falseIso}.

\medskip
\subsection{Delay and False Alarm Analysis}
The approximate delay and false alarm performance of the stopping rule $\tauhb$ is given in the theorem below. 
\begin{theorem}\label{thm:PerfTauHB}
Fix any $\epsilon_u > 0$ and $\epsilon_v > 0$. If the densities \eqref{eq:PDF_V_nok} and \eqref{eq:globalgU} are the true distributions of the samples then we have the following:
\begin{enumerate}
\item For the stopping rule $\tauhb$, the supremum in \eqref{prob:Pollak} is
achieved at $\gamma=1$, i.e.,
\[
\sup_{\gamma\geq 1} \; \Expect_\gamma[\tauhb-\gamma| \tauhb \geq \gamma] = \Expect_1[\tauhb-1].
\]
\item Setting $A_u=A_v=\log \beta$ ensures that as $\beta \to \infty$,
\[
\Expect_\infty[\tauhb] \geq \beta (1+o(1)).
\]
\item For each possible true post-change parameters $J$ and $\theta$, with $|J-1|\geq \epsilon_v$, and $|\theta - 1| \geq \epsilon_u$, as $\beta \to \infty$ 
\begin{equation}
\begin{split}
\Expect_1[\tauhb] &\leq  \left( \frac{\log \beta}{I(J)}  + \frac{\log \beta}{I(\theta)} \right)(1+o(1)),
\end{split}
\end{equation}
where $I(J)$ is the Kullback-Leibler divergence between the densities $f_V(\cdot; J)$ and $f_V(\cdot; 1)$, and 
$I(\theta)$ is the Kullback-Leibler divergence between the densities $g(\cdot; \theta)$ and $g(\cdot; 1)$.
\end{enumerate}
\end{theorem}
\begin{IEEEproof}
The proof follows from standard QCD arguments in the literature (see \cite{veer-bane-elsevierbook-2013}) and because
\[
\tauhb \geq \tauu,
\]
and 
\[
\tauhb \leq \tauu + \tauv. 
\]
\end{IEEEproof}

\section{Numerical Results}\label{sec:NumericalResults}
\begin{figure}[htb]
\center 
\includegraphics[width=9cm, height=6cm]{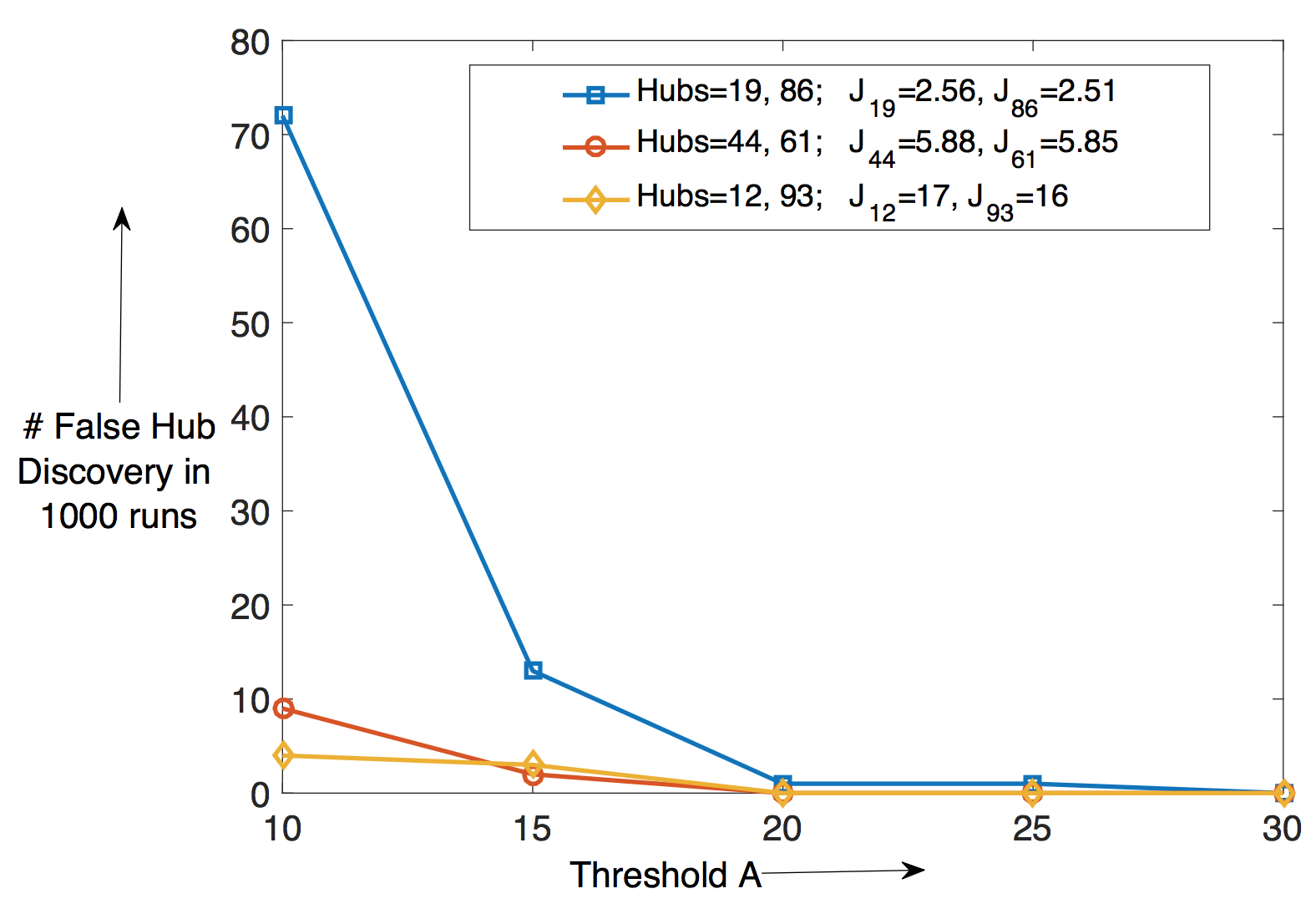}
\caption{False Isolation performance for three scenarios in \eqref{eq:FaultIsoCases}}
\label{fig:falseIso}
\end{figure}

In Fig.~\ref{fig:falseIso} we report the hub discovery performance of the proposed algorithm. 
We set $\mathbf{\Sigma}_0 = I_p$, the identity matrix. 
We assume that the post-change covariance matrix $\mathbf{\Sigma}$ is a row-sparse matrix of degree $j$,
obtained as follows. A $p \times p$ sample from the Wishart distribution is generated and some of
the entries are forced to be zero in such a way that no row has more than $j$ non-zero elements.
Specifically, we retain the top left $j \times j$ block of the matrix, and for each row $k$, $j+1 \leq k \leq (p+j)/2$,
all but the diagonal and the $(p+j+1-k)$th element is forced to zero. Each time an entry $(k,i)$ is set to zero,
the entry $(i,k)$ is also set to zero, to maintain symmetry. Finally, a scaled diagonal matrix is also added to $\mathbf{\Sigma}$
to restore its positive definiteness. We set $n=10$, $p=100$, and $j=5$.

To implement $\tauv$ we have chosen $\epsilon_v=1$, 
and we use the the maximum likelihood estimator for $J_k$ which, 
as a function of $m$ samples $(V_k(1), \cdots, V_k(m))$ from $f_V(\cdot, J)$, is given by
\begin{equation}\label{eq:MLEst_J}
\hat{J}(V_k(1), \cdots, V_k(m)) = \frac{1}{(p-1)\frac{1}{m} \sum_{i=1}^m P_0(V_k(i))}. 
\end{equation}
Specifically, 
\begin{equation}
\begin{split}
\argmax_{J: J \geq 2} \;  \log \sum_{i=\ell}^m & \frac{f_V(V_k(i);J)}{f_V(V_k(i); 1)} \\
&= \max\{2, \hat{J}(V_k(\ell), \cdots, V_k(m))\}.
\end{split}
\end{equation}
Similar techniques were used to implement $\tauu$; see \cite{bane-firo-hero-isit15}.

We consider three scenarios, each with two hubs:
\begin{equation}\label{eq:FaultIsoCases}
\begin{split}
1\;  \mathsf{H} &= \{19, 86\} \mbox{ with } J_{19} = 2.56 \mbox{ and } J_{86}=2.51.\\
2.\;  \mathsf{H} &= \{44, 61\} \mbox{ with } J_{44} = 5.88 \mbox{ and }  J_{61}=5.85.\\
3.\; \mathsf{H} &= \{12, 93\} \mbox{ with } J_{12} = 17 \mbox{ and }  J_{93}=16.
\end{split}
\end{equation}
The parameter $J_k$s are estimated using $1000$ samples. For each case, the total number of false 
isolations is calculated using $1000$ sample paths by assuming $\gamma=1$. The parameter $q$ was chosen to be $10$. 
The definition \eqref{eq:realDefHubDiscovery} was used as a criterion calculating false isolation or false hub discovery. The results for all the three cases is plotted in Fig.~\ref{fig:falseIso} as a function of the 
threshold $A=A_v=A_u$. As shown in Fig.~\ref{fig:falseIso}, the number of false discoveries goes to zero as the threshold is increased. Thus, empirical evidence suggests that
the test is consistent.

\begin{figure}[htb]
\center 
\includegraphics[width=9cm, height=6cm]{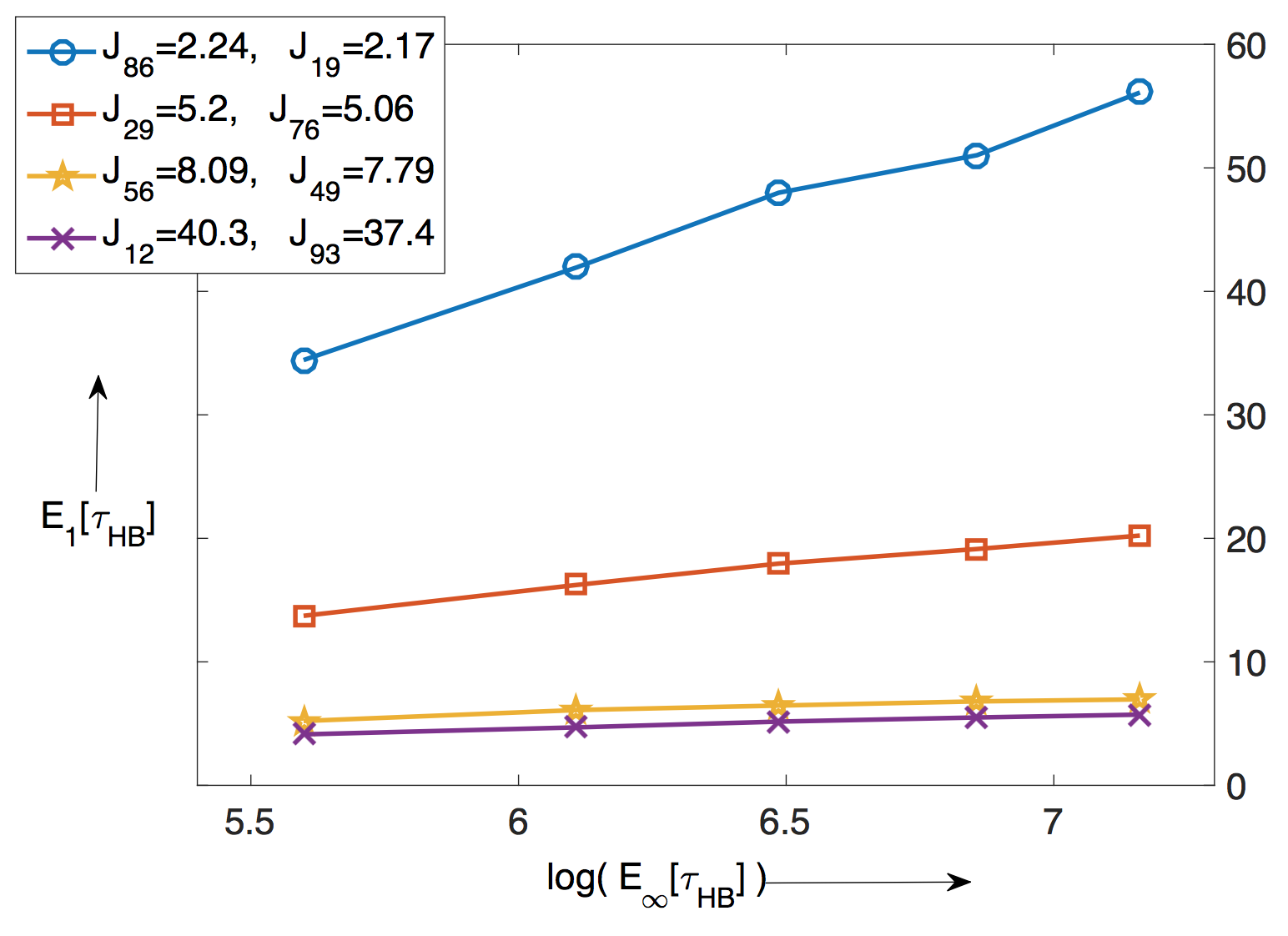}
\caption{MFA-Delay performance for the stopping rule $\tauhb$ for various choices of parameter $\{J_k\}$. 
The subscripts of $J$ are the hub indices.}
\label{fig:MFADelay}
\end{figure}

\medskip
In Fig.~\ref{fig:MFADelay} we plot the delay ($\Expect_1[\tauhb]$) vs the log of mean time to false alarm ($\log \Expect_\infty[\tauhb]$)
for various cases specified in the figure. The values in the figure are obtained
by choosing different values of the threshold $A$ and estimating the delay by choosing the change point $\gamma=1$
and simulating the test for $500$ sample paths.
The mean time to false alarm values are estimated by simulating the test for $1500$ sample paths.

\section{Conclusions and Future Work}
We have introduced the quickest hub detection and localization problem for high dimensional correlation graphs. Local and global summary statistics were proposed based on purely high dimensional asymptotic theory  for the presence of hub nodes, defined as vertices having either high degree or large correlations. The asymptotic densities of these summary statistics were shown to be in the one dimensional exponential family of densities. These densities were used  to define a hybrid local and global quickest detection test for hubs and the asymptotic performance of this test was evaluated in terms of mean time to detect a change and mean time to false alarm.

\section{Acknowledgments}
This work was partially supported by the Consortium for Verification Technology under Department of Energy National Nuclear Security Administration award number DOE-NA0002534.


\balance
{
\bibliographystyle{ieeetr}
\bibliography{QCD_verSubmitted}
}

\end{document}